\documentclass{aims}
\usepackage{amsmath}
\usepackage{paralist}
\usepackage{graphics} 
\usepackage{epsfig} 
\usepackage{graphicx}  \usepackage{epstopdf}
\usepackage[colorlinks=true]{hyperref}
\hypersetup{urlcolor=blue, citecolor=red}

\newtheorem{theorem}{Theorem}[section]

\newtheorem{lemma}[theorem]{Lemma}
\newtheorem{proposition}{Proposition}[section]

\theoremstyle{definition}
\newtheorem{definition}[theorem]{Definition}

\title[Uniqueness and stability of a fluid-like driven system] 
{Uniqueness and stability with respect to parameters of solutions to a fluid-like driven system for active-passive pedestrian dynamics}

\author[T. K. Thoa Thieu and Matteo Colangeli and Adrian Muntean]{}


\email{thoa.thieu@kau.se}
\email{matteo.colangeli1@univaq.it }
\email{adrian.muntean@kau.se}

\thanks{$^*$ Corresponding author: thikimthoa.thieu@kau.se}

\begin{document}
	\maketitle
	
	\centerline{\scshape T. K. Thoa Thieu$^*$}
	\medskip
	{\footnotesize
		\centerline{Department of Mathematics, Gran Sasso Science Institute,  }
		\centerline{Viale Francesco Crispi 7, L'Aquila 67100, Italy}
		\centerline{Department of Mathematics and Computer Science, Karlstad University,}
		\centerline{Universitetsgatan 2, Karlstad, Sweden}
	} 
	
	\medskip
	
	\centerline{\scshape Matteo Colangeli}
	\medskip
	{\footnotesize
		\centerline{ Dipartimento di Ingegneria e Scienze dell'Informazione e Matematica, Universit\`{a} degli Studi dell'Aquila,}
		\centerline{Via Vetoio, L'Aquila 67100, Italy}
	}
	
	\medskip
	
	\centerline{\scshape Adrian Muntean}
	\medskip
	{\footnotesize
		\centerline{Department of Mathematics and Computer Science, Karlstad University,}
		\centerline{Universitetsgatan 2, Karlstad, Sweden}
	}
	\bigskip
	

	\begin{abstract}
		We study a system of parabolic equations consisting of a double nonlinear parabolic equations of Forchheimer type coupled with a semilinear parabolic equations. The system describes a fluid-like driven system for active-passive pedestrian dynamics. The structure of the nonlinearity of the coupling allows us to prove the uniqueness of solutions. We provide also stability estimates of solutions with respect to selected parameters.
	\end{abstract}
	
	\section{Introduction}\label{Sec:intro}	
	A recent result on the weak solvability of a mixed fluid-like driven system for active--passive has been reported in \cite{Thieu2019}, where the authors provided the existence of solutions to the problem \eqref{main_eq} by using a Schauder's fixed point argument. This type of mixed pedestrian dynamics is originally proposed in \cite{Richardson2019} by considering their evacuation dynamics in a complex geometry in the presence of a fire as well as of a slowly spreading smoke curtain. From a stochastic processes perspective, various lattice gas models for active-passive pedestrian dynamics have been explored in \cite{Cirillo2019, Colangeli2019}. Within the present framework, our model is embedded in a continuum scale and resembles the structure of Forchheimer flows in porous media \cite{Bear1972}. The aim of this paper is to complete the proof of the well-posedness of the system \eqref{main_eq} by showing the uniqueness and stability of solutions with respect to parameters.  The nonlinear structure in the transport term where the Forchheimer polynomial appears, allows us to establish the wanted uniqueness and stability estimates. This work focuses on the structural stability of solution with respect to initial and boundary data, nonlinear coupling coefficient, and to the difussion coeficient from the semi-linear equation. We also provide a stability-like estimate for the gradient of the solutions. To obtain such control on the gradients, the structure of the system has an important degenerate monotonicity property\footnote{This terminalogy is taken from \cite{Hoang2011} and refers to the assumption $(\text{A}_5)$.}, which allows us to compare the difference between pairs of parametrized solutions and their gradients (see in  \cite{Hoang2011}). 
    
    A number of relevant results are available on structural stability topics. In particular, standard nonlinear energy stability results have been presented in \cite{Straughan2014} for convection problems, where the author dealt with an integral inequality technique referred to as the energy method. The structral stability of solutions to generalized Forchheimer equations (introduced in \cite{Celik2016}) has been provided in \cite{Aulisa2009}, where the authors investigated the uniqueness, the Lyapunov asymptotic stability together with the large time behavior features of the corresponding initial boundary value problems. A structural stability with respect to boundary data and the coefficients of Forchheimer is considered in \cite{Hoang2011}. In \cite{Wijaya2019}, a stability estimate is introduced by considering a nonlinear drag force term corresponding to the Forchheimer term in a Navier–Stokes type model of flow in non-homogeneous porous media. Such investigations on stability estimates not only contribute to the understanding of the well-posedness of model equations, but also can point out inherent delimitations of the parameters regions outside which it makes no sense to search for solutions, see e.g. \cite{Vromans2019}.
    
    This paper is organized as follows. In Section \ref{Sec:Pre}, preliminaries and assumptions are provided. Then, we recall available energy estimates in Section \ref{Sec:Energy}. In Section \ref{uniqueness}, we show the proof of the uniqueness of solutions to our system. Finally, the target stability estimates are established in Section \ref{stability} and Section \ref{stability-full}.

	
	\section{Setting of the model equations}
	Let a bounded set $\Omega \neq \emptyset$, $\Omega \subseteq \mathbb{R}^2$ has $C^1$-boundary\footnote{This boundary information is to guarantee the trace's inequality \eqref{trace-theo} (e.g. \cite{Celik2016}).} $\partial \Omega$ such that $\partial \Omega=\Gamma^N \cup \Gamma^R$, $\Gamma^N \cap \Gamma^R= \emptyset$ with $\mathcal{H}(\Gamma^N) \neq \emptyset$ and $\mathcal{H}(\Gamma^R)\neq \emptyset$, where $\mathcal{H}$ denotes the surface measure on $\Gamma^N, \Gamma^R$ and take $S=(0,T)$. We shall consider the following equations, where the pair of velocities is $(u=u(t,x),v=v(t,x))$ such that the mappings $u: S \times \Omega \longrightarrow \mathbb{R}^2$ and $v: S\times\Omega \longrightarrow \mathbb{R}^2$ satisfy
	\begin{align}\label{main_eq}
	\begin{cases}
	\partial_t (u^\lambda) + \mathrm{div}(-K_1(|\nabla u|)\nabla u) = -b(u-v) \ \text{ in } S \times \Omega,\\
	\partial_t v - K_2\Delta v =b(u-v) \ \text{ in } S \times \Omega, \\
	-K_1(|\nabla u|)\nabla u \cdot {\bf n} = \varphi u^\lambda \ \text{ at } S \times\Gamma^R,\\
	-K_1(|\nabla u|)\nabla u\cdot {\bf n} = 0 \ \text{ at } S\times\Gamma^N ,\\
	-K_2\nabla v \cdot {\bf n} = 0 \ \text{ at }  S\times \partial\Omega, \\
	u(t=0,x) = u_0(x) \ \text{ for } x \in \bar{\Omega}, \\
	v(t=0,x) = v_0(x) \ \text{ for } x \in \bar{\Omega}.
	\end{cases}
	\end{align}
	In \eqref{main_eq}, $K_2 >0$ and function $K_1$ stems from the derivation of a nonlinear version of the Darcy equation defined via a generalized polynomial with non-negative coefficients (e.g. \cite{Hoang2011}, \cite{Aulisa2009}, \cite{Celik2016}). The structure of $K_1$ in the first equation of \eqref{main_eq} will be described in Section $2$. In addition, $\lambda \in (0,1]$ is a fixed number and $b(\cdot)$ is a sink/source term. The system \eqref{main_eq} has been studied in \cite{Thieu2019}.
	
	\section{Preliminaries and Assumptions}\label{Sec:Pre}
	We list in this section a couple of preliminary results (mostly inequalities and compactness results) as well as our assumption on data and parameters.
	\begin{lemma} Let $x, y \geq 0$. Then the following elementary inequalities hold:
		\begin{equation}\label{pre-ine-1}
		(x+y)^p \leq 2^p(x^p+y^p) \text{ for all } p >0, 
		\end{equation}
		\begin{equation}\label{pre-ine-2}
		(x+y)^p \leq x^p+y^p \text{ for all } 0 <p \leq 1,
		\end{equation}
		\begin{equation}\label{pre-ine-3}
		(x+y)^p \leq 2^{p-1}(x^p+y^p) \text{ for all } p \geq 1,
		\end{equation}
		\begin{equation}\label{pre-ine-4}
		x^{\beta} \leq x^{\alpha} + x^{\gamma} \text{ for all } 0 \leq \alpha \leq \beta \leq \gamma,
		\end{equation}
		\begin{equation}\label{pre-ine-5}
		x^\beta \leq 1+x^{\gamma} \text{ for all } 0 \leq \beta \leq \gamma.
		\end{equation}
	\end{lemma}
	The proof is elementary and we omit it from here.
	\begin{lemma}[\bf Trace Lemma]\label{trace-lm}
		Let $\lambda \in (0,1]$, $\delta =1-\lambda$, $a = \frac{\alpha_N}{\alpha_N+1} \in (0,1)$, $a>\delta$, $\alpha \geq 2-\delta$, $\alpha \leq 2$, $\mu_0=\frac{a-\delta}{1-a}$, $\alpha_{\star}=\frac{n(a-\delta)}{2-a}$ and 
		\begin{align}
		\theta = \theta_\alpha := \frac{1}{(1-a)(\alpha/\alpha_{\star}-1)} \in (0,1).
		\end{align}
		Then it exists $C>0$ such that the following estimate holds
		\begin{align}\label{trace-theo}
		\int_{\Gamma^R}|u|^\alpha d\sigma \leq 2\varepsilon \int_{\Omega} |u|^{\alpha+\delta-2}|\nabla u|^{2-a}dx + C\|u\|_{L^{\alpha}(\Omega)}^\alpha \nonumber\\+ C\varepsilon^{-\frac{1}{1-a}}\|u\|_{L^{\alpha}(\Omega)}^{\alpha+\mu_0} + C\varepsilon^{-\mu_2}\|u\|_{L^\alpha(\Omega)}^{\alpha+\mu_1},
		\end{align}
		where
		\begin{align}
		\mu_1&=\mu_{1,\alpha} := \frac{\mu_0(1+\theta(1-a))}{1-\theta},\\
		\mu_2&=\mu_{2,\alpha}:=\frac{1}{1-a}+\frac{\theta(2-a)}{(1-\theta)(1-a)}.
		\end{align}
	\end{lemma}
	For the proof of Lemma $2.2$, see Lemma $2.2$  in \cite{Celik2016}.
	%
	\subsection{Structure of $K_1$}
	In this section, we recall the definitions on the constructions based on the nonlinear Darcy equation and its monotonicity properties as they have been presented in \cite{Aulisa2009}.
	First of all, we introduce the function $K_1: \mathbb{R}^{+}\longrightarrow \mathbb{R}^{+}$ defined for $\xi \geq 0$ by $K_1(\xi)=\frac{1}{g(s(\xi))}$ which is supported to be the unique non-negative solution of the equation $sg(s) = \xi$, where $g :\mathbb{R}^{+} \to \mathbb{R}^{+}$ is a polynomial with positive coefficients defined by 
	\begin{align}\label{g}
	g(s) = a_0s^{\alpha_0} + a_1s^{\alpha_1} + \ldots + a_Ns^{\alpha_N} \ \text{ for } s \geq 0,
	\end{align}
	where $\alpha_k \in \mathbb{R}_{+}$ with $k\in \{0,\ldots,N\}$.\\
	The function $g$ is taken to be independent of the spatial variable. Thus, we may have
	\begin{align}\label{nonlinearDarcy}
	G(|v|)=g(|v|)|v|=|\nabla p|,
	\end{align}
	where $G(s)=sg(s)$ for $s \geq 0$. From now on we use the following notation for the function $G$ and its inverse $G^{-1}$, namely, $G(s) = sg(s)=\xi$ and $s=G^{-1}(\xi)$. To be successful with the analysis to follow, we impose the following condition on the polynomial $g$, referred to as $(G)$.

	\textit{$(G_1)$ } $g \in C([0,\infty)) \cap C^1((0,\infty))$ such that
	$$g(0) > 0 \text{ and } g'(s) \geq 0 \text{ for all } s \geq 0.$$
	
	\textit{$(G_2)$ } It exists $\theta >0 $ with $g \in C([0,\infty)) \cap C^1((0,\infty))$ such that
	\begin{align}
	g(s) \geq \theta sg'(s) \text{ for all } s>0.
	\end{align}
	
	To be able to ensure the uniqueness of solution to the system \eqref{main_eq}, we use the monotonicity properties of the function $F: \mathbb{R}^d \longrightarrow \mathbb{R}^d$ such that $F(y)= K_1(|y|)y$. This is related to the nonlinear Darcy structure \eqref{nonlinearDarcy}. Furthermore, we recall the following basic essential ingredients:
	\begin{definition}
		Let $F: \mathbb{R}^d \longrightarrow \mathbb{R}^d$ be a given mapping.
		\begin{itemize}
			\item $F$ is monotone if 
			\begin{align}
			(F(y') - F(y))\cdot (y'-y) \geq 0 \text{ for all } y',y \in \mathbb{R}^d.
			\end{align}
			\item $F$ is strictly monotone if there is $c>0$, such that
			\begin{align}
			(F(y') - F(y))\cdot (y'-y) \geq c|y'-y|^2 \text{ for all } y',y \in \mathbb{R}^d.
			\end{align}
			\item $F$ is strictly monotone on bounded sets if for any $R>0$, there is a positive number $c_R >0$, such that
			\begin{align}
			(F(y')-F(y))\cdot (y'-y) \geq c_R|y'-y|^2 \text{ for all } |y'| \leq R, |y| \leq R.
			\end{align}
		\end{itemize}
		See Definition III.3 in \cite{Aulisa2009} for more details.
	\end{definition}
	
	We introduce a useful formulation by defining the following function $\Phi: \mathbb{R}^d \times \mathbb{R}^d \longrightarrow \mathbb{R}$ as follow
	\begin{align}\label{Phi}
	\Phi(y,y') = (K_1(|y'|)y' - K_1(|y|)y)\cdot (y'-y) \text{ for } y,y' \in \mathbb{R}^d.
	\end{align} 
	\begin{proposition}\label{prop-phi}
		Let $g$ satisfy $(G_1)$.
		Then $F(y) = K_1(|y|)y$ is monotone, hence $\Phi(y,y') \geq 0$ for all $y,y' \in \mathbb{R}^d$, where $\Phi$ is defined as in \eqref{Phi}. 
	\end{proposition}
	For the proof of Proposition \ref{prop-phi}, see Proposition III.4 in \cite{Aulisa2009}.
	\begin{lemma}\label{decrease}
		Let $g$ satisfies $(G_1)$. The function $K_1(\cdot)=K_{1g}(\cdot)=\frac{1}{g(s(\cdot))}$, is well defined, belongs to $C^1([0,\infty))$, and is decreasing. Moreover, for any $\xi \geq 0$, let $s=G^{-1}(\xi)$, then one has
		\begin{align}
		K_1'(\xi) = -K_1(\xi)\frac{g'(s)}{\xi g'(s) + g^2(s)} \leq 0.
		\end{align}
	\end{lemma}
	For the proof of Lemma \ref{decrease}, see Lemma III.2 in \cite{Aulisa2009}.
	\begin{proposition}\label{prop-phi2}
		Let $g$ satisfies $(G_1)$ and $(G_2)$. Then $F(y)=K_1(|y|)y$ is strictly monotone on bounded sets. More precisely,
		\begin{align}\label{monotonicity}
		\Phi(y,y') \geq \frac{\lambda}{\lambda+1}K_1(\max\{|y|,|y'|\})|y'-y|^2 \text{ for all } y,y' \in \mathbb{R}^d.
		\end{align}
	\end{proposition}
	For the proof \ref{prop-phi2}, see Proposition III.6 in \cite{Aulisa2009}.
	\subsection{Assumptions}
	We make the following choices on the structure of the involved nonlinearities.
	\begin{itemize}
		\item[($\text{A}_1$)] The structure of $K_1(\xi)$ has the following properties hold $K_1: [0,\infty) \longrightarrow (0,\frac{1}{a_0}]$ such that $K_1$ is decreasing and
		\begin{align}\label{k1-eqn1}
		\frac{d_1}{(1+\xi)^a} \leq K_1(\xi) \leq \frac{d_2}{(1+\xi)^a};
		\end{align}
		\begin{align}\label{k1-eqn2}
		d_3(\xi^{2-a}-1)\leq K_1(\xi)\xi^2 \leq d_2\xi^{2-a} \text{ for all } \xi \in [0,\infty).
		\end{align}
		In \eqref{k1-eqn1}, $d_1, d_2, d_3$ are strictly positive constants depending on $g(s)$ and $a \in (0,1)$. 
		\item[($\text{A}_2$)] The function $b: \mathbb{R} \longrightarrow \mathbb{R}$ satisfies the following structural condition: it exists $\hat{c}>0$ such that $b(z) \leq \hat{c}|z|^{\sigma}$, with $\sigma \in (0,1)$. 
		\item[($\text{A}_3$)] The source term $b: \mathbb{R} \longrightarrow \mathbb{R}$ is globally Lipschitz continuous.
		\item[($\text{A}_4$)]  The boundary data satisfies $\varphi \in L^{\infty}(\Gamma^N)$.
		\item[($\text{A}_5$)] $K_1$ satisfies the degeneracy type
		\begin{align*}
		C_1(1+\xi)^{-a} \leq K_1(\xi) \leq C_2(1+\xi)^{-a},
		\end{align*}
		where $C_1$ and $C_2$ are constants.
	\end{itemize}
	Assumptions $(\text{A}_1)$-$(\text{A}_5)$ are all technical.The choice of $(\text{A}_1)$ was inspired by Theorem III.10 in \cite{Aulisa2009} and the choice of $(\text{A}_5)$ corresponds to the setting in \cite{Hoang2011}.
	
	We recall from \cite{Thieu2019} the following concept of solution to \eqref{main_eq} fitting to the case $\alpha \in [1+\lambda,2]$.
	\begin{definition}\label{dn-weak}
		Find $$(u,v) \in L^\alpha(S;L^\alpha(\Omega))\cap L^{2-a}(S;W^{1,2-a}(\Omega)) \times L^2(S;W^{1,2}(\Omega))$$
		satisfying the identities
		\begin{align}\label{weak1}
		\int_{\Omega}\partial_t(u^{\lambda})\psi dx + \int_{\Omega} K_1(|\nabla u|)\nabla u \nabla \psi dx + \int_{\Gamma^R}\varphi u^{\lambda}\psi d\gamma = -\int_{\Omega}b(u-v)\psi dx
		\end{align}
		and
		\begin{align}\label{weak2}
		\int_{\Omega}\partial_t v\phi dx + \int_{\Omega}K_2 \nabla v \nabla \phi dx = \int_{\Omega}b(u-v)\phi dx
		\end{align}
		for a.e $t \in S$ and for all $(\psi, \phi) \in L^\alpha(\Omega) \times  W^{1,2}(\Omega)$ with the initial data $(u_0,v_0) \in L^\alpha(\Omega) \times L^2(\Omega)$.
	\end{definition}
	This weak formulation has been presented in \cite{Thieu2019}.
	\begin{theorem}\label{existence-theo}
		Assume that $(A_1)$ and $(A_2)$ hold. Let $\lambda \in (0,1]$, $\delta =1-\lambda$, $a = \frac{\alpha_N}{\alpha_N+1} \in (0,1)$, $a > \delta$, $\alpha \geq 2-\delta$, $\alpha \leq 2$, $\sigma\leq \frac{\alpha}{2}$, $\sigma \in (0,1)$ and $u_0 \in L^{\alpha}(\Omega)$, $v_0 \in L^2(\Omega)$. Then the problem \eqref{main_eq} has at least a weak solution\\ $(u,v) \in L^\alpha(S;L^\alpha(\Omega))\cap L^{2-a}(S;W^{1,2-a}(\Omega)) \times L^2(S;W^{1,2}(\Omega))$ in the sense of Definition \ref{dn-weak}. 
	\end{theorem}
	For the proof of this result, see Theorem $3.4$ in \cite{Thieu2019}.
	\subsection{Statement of the main results}
	The main results of this paper are stated in Theorem \ref{uniqueness-theo}, Theorem \ref{stability-theo} and Theorem \ref{stabilityfull-theo}. They correspond to the case $\lambda=1, \delta=0$.
	\begin{theorem}\label{uniqueness-theo}
		Assume that $(\text{A}_1)$-$(\text{A}_4)$ hold. Let $a = \frac{\alpha_N}{\alpha_N+1} \in (0,1)$, $a > \delta$, $\alpha \geq 2-\delta$, $\alpha \leq 2$ and $\sigma\leq \frac{\alpha}{2}$, $\sigma \in (0,1)$. Then, the problem \eqref{main_eq} admits at most a weak solution in the sense of Definition \ref{dn-weak}.
	\end{theorem}
	We look for the case when the coupling is linear, i.e $b: \mathbb{R} \longrightarrow \mathbb{R}$ is a given function such that $b(s)=rB(s)$, where $r\in (0,\infty)$. Here, $B$ is taken such that $(\text{A}_3)$ and $(\text{A}_4)$ are satisfied. We call $S_1=(0,T_1)$, $S_2=(0,T_2)$ and $
	\mathcal{S}=(0,\min\{T_1,T_2\}) = (0,\tau)$. 
	Let us define a triplet \\ $(u_i, v_i, \mathcal{D}_i)$, where $(u_i,v_i) \in \left(L^\alpha(S;L^\alpha(\Omega))\cap L^{2-a}(S;W^{1,2-a}(\Omega))\right) \times L^2(S;W^{1,2}(\Omega))$ and $\mathcal{D}_i= (D_i, r_i, u_{0i}, v_{0i}) \in (0,\infty) \times (0,\infty) \times L^\alpha(\Omega) \times L^2(\Omega)$. To avoid the use of multiple indices, we denote $D:= K_2$, where $K_2>0$ is entering \eqref{main_eq}. We give stability estimates of the solutions with respect to initial and boundary data, nonlinear coupling coefficient $r$ and the diffusion coefficient $D$. 
	\begin{theorem}\label{stability-theo}
		Assume that $(\text{A}_1)$-$(\text{A}_4)$ hold. For $i\in\{1,2\}$, $(D_i, r_i, u_{0i}, v_{0i})$ belong to a fixed compact subset $K \subset (0,\infty) \times (0,\infty) \times L^\alpha(\Omega) \times L^2(\Omega)$,  $\lambda=1$, $\bar{r} \geq |r_1-r_2|> 0$. Let $(u_i, v_i)$ be weak solutions to \eqref{main_eq} corresponding to the choices of data $(D_i, \varphi_i, r_i, u_{0i}, v_{0i}), i\in\{1,2\}$. Then, the following stability estimate holds
		\begin{align}\label{stabi-main1}
		\|u_1-u_2\|_{L^\alpha(\Omega)}^\alpha + \|v_1-v_2\|_{L^2(\Omega)}^2 \leq e^{C(\alpha,\lambda,\hat{c},\bar{r})|r_1-r_2|t}\Big[ \|u_{01} - u_{02}\|_{L^\alpha(\Omega)}^\alpha \nonumber\\+ \|v_{01} - v_{02}\|_{L^2(\Omega)}^2 + Ct(|D_1-D_2|+ |r_1-r_2| - \|\varphi_1-\varphi_2\|_{L^\infty(\Gamma^R)}^2 )\Big],
		\end{align}
		for $t\in \mathcal{S}$.
	\end{theorem}
	\begin{theorem}\label{stabilityfull-theo}
		Assume that $(\text{A}_1)$-$(\text{A}_5)$ hold. For $i\in \{1,2\}$, let $(D_i, r_i, u_{0i}, v_{0i})$ belong to a fixed compact subset $K \subset (0,\infty) \times (0,\infty) \times L^\alpha(\Omega) \times L^2(\Omega)$, $\bar{r} \geq |r_1-r_2|> 0$. Let $(u_i, v_i)$ be weak solutions to \eqref{main_eq} corresponding to the choices of data $(D_i, r_i, u_{0i}, v_{0i})$. Then, the following estimate holds
		\begin{align}\label{stabi-main2}
		\int_{\Omega}|\nabla u_1-\nabla u_2|^{2-a}dx + \int_{\Omega}|\nabla v_1-\nabla v_2|^2dx \leq C+ C(|D_1-D_2| -\nonumber\\ \|\varphi_1-\varphi_2\|_{L^\infty(\Gamma^R)}^2)  + \left(\left(t + \frac{1}{2}\right) + |r_1-r_2|\frac{\hat{c}}{2\bar{r}}\right)e^{C(\alpha,\lambda,\hat{c},\bar{r})|r_1-r_2|t}\Big(\|u_{01} - u_{02}\|_{L^2(\Omega)}^2 \nonumber\\ + \|v_{01}-v_{02}\|_{L^2(\Omega)}^2+ Ct\left(|D_1-D_2|+|r_1-r_2| - \|\varphi_1-\varphi_2\|_{L^\infty(\Gamma^R)}^2\right)\Big).
		\end{align}
		for $t \in \mathcal{S}$.
	\end{theorem}
	The proofs of Theorem \ref{uniqueness}, Theorem \ref{stability-theo} and Theorem \ref{stabilityfull-theo} are given in  Section \ref{uniqueness}, Section \ref{stability} and Section \ref{stability-full}, respectively.
	
	Note that our stability estimate for the gradient of the solutions are not optimal. The bound \eqref{stabi-main2} can be eventually improved by studying the structural stability with respect to the Forchheimer polynomial $K_1(\cdot)$. We refer the reader to \cite{Hoang2011}, where this case has been studied for a simpler setting.
	\section{Energy estimates}\label{Sec:Energy}
	In this section, we recall the energy estimates available for the problem \eqref{main_eq}. In particular, Proposition \ref{prop1} contains $L^\alpha-L^2$ estimates, while gradient and time derivative estimates are reported in Proposition \ref{prop2}.
	\begin{proposition}\label{prop1}
		Assume that $(\text{A}_1)$-$(\text{A}_5)$ hold and let $\lambda \in (0,1]$, $\delta =1-\lambda$, $a = \frac{\alpha_N}{\alpha_N+1} \in (0,1)$, $a > \delta$, $\alpha \geq 2-\delta$, $\alpha \leq 2$, $\sigma\leq \frac{\alpha}{2}$, $\sigma \in (0,1)$ and $u_0 \in L^{\alpha}(\Omega)$, $v_0 \in L^2(\Omega)$. Then, for any $t\in S$, the following estimates hold
		\begin{align}\label{main_0}
		\frac{d}{dt}\int_{\Omega}|u|^{\alpha}dx + \int_{\Omega}|\nabla u|^{2-a}|u|^{\alpha + \delta -2}dx \leq C_1+\Big(\frac{3}{2C_2}\hat{c}+\nonumber\\ + \frac{d_3(\alpha - \lambda)}{C_2}\Big)\|u\|_{L^{\alpha}(\Omega)}^{\alpha} + \frac{\hat{c}}{2C_2}\|v\|_{L^2(\Omega)}^2.
		\end{align}
		\begin{align}\label{ine_main1}
		\int_{\Omega}|u|^{\alpha}dx + \int_{\Omega}v^2dx \leq e^{C_3t}\left(1+\|u_0\|_{L^\alpha(\Omega)}^{\alpha} + \|v_0\|_{L^2(\Omega)}^2\right),
		\end{align}
		\begin{align}\label{ine_main2}
		\int_0^T\int_{\Omega}|u|^{\alpha+\delta-2}|\nabla u|^{2-a}dxdt + \int_{0}^{T}\int_{\Omega}|\nabla v|^2dxdt \leq C_5 +\nonumber\\+ C_6\left(\|u_0\|_{L^\alpha(S; L^\alpha(\Omega))}^{\alpha} + \|v_0\|_{L^2(S; L^2(\Omega))}^2 \right),
		\end{align}
		where $C_1:= \frac{d_3(\alpha-\lambda)+\hat{c}}{C_5}|\Omega|$, $C_2:=\min\left\{\frac{\lambda}{\alpha},d_3(\alpha-\lambda)\right\}$, $C_3:=\max\left\{\frac{5}{2}\underline{\tilde{c}}\hat{c}|\Omega|,2\underline{\tilde{c}}\hat{c}\right\}$, $C_4:= \min\left\{\alpha-\lambda,K_2\right\}$,  $C_5:=\frac{5T}{2C_4}\hat{c}|\Omega| + \frac{2T\hat{c}e^{C_3t}}{C_4}$,  and   $C_6:=\frac{2\hat{c}e^{C_3t}}{C_2}$ with $\underline{\tilde{c}}:=\frac{1}{\min\left\{\frac{\lambda}{\alpha},\frac{1}{2}\right\}}$ and $\hat{c}$ is as in $(\text{A}_2)$, respectively.
	\end{proposition}
	For the proof of this result, see Proposition $4.1$ in \cite{Thieu2019}.
	We consider the following function $H: \mathbb{R}_{+} \longrightarrow \mathbb{R}_{+}$ given by
	\begin{align}
	H(\xi) = \int_{0}^{\xi^2}K_1(\sqrt{s})ds  \text{ for } \xi \in \mathbb{R}_{+}.
	\end{align}
	We admit a structural inequality between  $H(\xi)$ and $K_1(\xi)\xi^2$ of the form:
	\begin{align}\label{k1-ineq2}
	K_1(\xi)\xi^2 \leq H(\xi) \leq 2K_1(\xi)\xi^2 \text{ for all } \xi \in \mathbb{R}_{+}.
	\end{align}
	By combining \eqref{k1-eqn1} and \eqref{k1-ineq2}, we deduce also that
	\begin{align}\label{H-ine1}
	d_3(\xi^{2-a}-1) \leq H(\xi) \leq 2d_2\xi^{2-a} \text{ for all } \xi \in \mathbb{R}_{+}. 
	\end{align}
    In \eqref{H-ine1}, $d_2, d_3$ and $a$ are defined as in $(\text{A}_1)$.
	\begin{proposition}\label{prop2}
		Assume that ($\text{A}_1$) and ($\text{A}_2$) hold. Let $\lambda \in (0,1]$, $\delta =1-\lambda$, $a=\frac{\alpha_N}{\alpha_N+1}$, $a > \delta$, $\alpha 
		\geq 2-\delta$, $\alpha\leq 2$ and $\sigma \leq \frac{\alpha}{2}$. Furthermore, suppose that $\nabla u_0 \in L^{\alpha}(\Omega) \cap L^{2-a}(\Omega)$, $u_0 \in L^{\alpha}(\Omega)$, $v_0 \in H^1(\Omega)$ and $\varphi \in L^{\infty}(\Gamma^R)$. Then, for any $t\in S$, the following estimates hold
		\begin{align}\label{ine-main3}
		\int_{\Omega}|\nabla u|^{2-a}dx + \int_{\Omega}|\nabla v|^2dx \leq C( \hat{c}, \lambda,a)\Bigg[\Lambda(0) + \int_{0}^{t}(1+\|u\|_{L^{\alpha}(\Omega)}^{\alpha})^{\beta}ds \nonumber\\+ \int_{0}^{t}\|v\|_{L^2(\Omega)}^2ds + \int_{0}^{t}\int_{\Gamma^R}|\varphi_t|^{\frac{\alpha}{\alpha-\lambda-1}}d\sigma ds\Bigg] + \int_{\Omega}|\nabla v_0|^2dx\nonumber\\
		+ \frac{\hat{c}^2}{C_2}|\Omega|t + \frac{\hat{c}^2}{2C_2} e^{C_3t}\left(1+\|u_0\|_{L^\alpha(S,L^{\alpha}(\Omega))}^{\alpha} + \|v_0\|_{L^2(S,L^2(\Omega))}^2\right) .
		\end{align}
		\begin{align}\label{ine-main4}
		\int_{\Omega}|(u^{\lambda})_t|^2dx + \int_{\Omega}|v_t|^2dx \leq C(\hat{c}, \lambda,a)\Bigg[1 + (1+\|u\|_{L^{\alpha}(\Omega)}^{\alpha})^{\beta} + \|v\|_{L^2(\Omega)}^2 \nonumber\\ + \int_{\Gamma^R}|\varphi_t|^{\frac{\alpha}{\alpha-\lambda-1}}d\sigma\Bigg] + \frac{\hat{c}^2}{C_2}|\Omega| + \frac{\hat{c}^2}{2C_2}e^{C_3t}\left(1+\|u_0\|_{L^{\alpha}(\Omega)}^{\alpha} + \|v_0\|_{L^2(\Omega)}^2\right),
		\end{align}
		where $C( \hat{c}, \lambda,a)>0$ is a constant and $$\Lambda(0):= \frac{\lambda+1}{2}\int_{\Omega}H(|\nabla u_0|)dx  + \int_{\Omega}|u_0|^{\alpha}dx.$$
	\end{proposition}
	For the proof of this result, see Proposition $4.2$ in \cite{Thieu2019}.
	\section{Proof of Theorem \ref{uniqueness-theo}}\label{uniqueness}
	\begin{proof}
		To prove the uniqueness of solutions in the sense of Definition \ref{dn-weak}, we adapt the arguments by E. Aulisa et. al (cf. Section IV, \cite{Aulisa2009}) to our setting. Essentially, we are using the monotonicity properties of the term $K_1(y)y$ as stated in Proposition \ref{prop-phi} and Proposition \ref{prop-phi2}.
		
		Let $(u_i, v_i), i \in \{1,2\}$ be two arbitrary weak solutions to  problem \eqref{main_eq} in the sense of Definition \ref{dn-weak}, where the initial data is take $u_i(t=0,x)=u_{i0}(x)$ and $v_i(t=0,x)=v_{i0}(x)$ for all $x \in \bar{\Omega}$. We denote $w=u_1-u_2$ and $z=v_1-v_2$. If we substitute the pair $(w,z)$ into \eqref{weak1}-\eqref{weak2}, we obtain
		\begin{align}
		\int_{\Omega} \partial_t (u_1^\lambda - u_2^\lambda)\psi dx + \int_{\Omega}\partial_t z \phi dx + \int_{\Omega} \left(K_1(|\nabla u_1|)\nabla u_1 - K_1(|\nabla u_2|)\nabla u_2\right)\nabla \psi dx + \nonumber\\ +K_2\int_{\Omega} \nabla z\nabla \phi dx=-\int_{\Gamma^R} \varphi (u_1^\lambda - u_2^\lambda)\psi d\gamma - \int_{\Omega}(b(u_1-v_1) - b(u_2-v_2))(\psi - \phi)dx
		\end{align}
		Now, choosing the test function $$(\psi,\phi):= (|w|^{\alpha+\delta-1},z)\in \left(\left(L^\alpha(\Omega)\cap W^{1,2-a}(\Omega)\right) \times W^{1,2}(\Omega)\right)$$ leads to
		\begin{align}\label{ine-u1}
		&\frac{\lambda}{\alpha}\frac{d}{dt}\int_{\Omega}|w|^\alpha dx +\frac{1}{2}\frac{d}{dt}\int_{\Omega}z^2dx+\nonumber\\ &\int_{\Omega} \Big(K_1(|\nabla u_1|)\nabla u_1 - K_1(|\nabla u_2|)\nabla u_2\Big)\nabla w |w|^{\alpha+\delta-2} dx + K_2\int_{\Omega}|\nabla z|^2dx \nonumber\\& +\int_{\Gamma^R} \varphi (u_1^\lambda - u_2^\lambda)|w|^{\alpha+\delta-1} d\gamma =- \int_{\Omega}(b(u_1-v_1)-b(u_2-v_2))(|w|^{\alpha+\delta-1} - z)dx.
		\end{align}
		Using assumption $(\text{A}_3)$ to handle the right hand side of \eqref{ine-u1}, we have the following estimate 
		\begin{align}\label{ine-u2}
		&\frac{\lambda}{\alpha}\frac{d}{dt}\int_{\Omega}|w|^\alpha dx +\frac{1}{2}\frac{d}{dt}\int_{\Omega}z^2dx+  \int_{\Omega} \Phi(\nabla u_1,\nabla u_2) |w|^{\alpha+\delta-2} dx \nonumber\\  &+ K_2\int_{\Omega}|\nabla z|^2dx + \int_{\Gamma^R} |\varphi||u_1^\lambda - u_2^\lambda| |w|^{\alpha+\delta-1} d\gamma  \nonumber\\&\leq  \left|\int_{\Omega}(b(u_1-v_1)-b(u_2-v_2))(|w|^{\alpha+\delta-1} - z)dx\right| \nonumber\\
		&\leq  \left|\int_{\Omega}(|u_1-u_2|+|v_1-v_2|)(|w|^{\alpha+\delta-1} - z)dx\right| \nonumber\\
		&\leq  \left|\int_{\Omega}|w|^{\alpha+\delta}dx - \int_{\Omega}|w||z|dx - \int_{\Omega}|v_1-v_2||w|^{\alpha+ \delta-1}dx + \int_{\Omega}z^2dx\right| \nonumber\\
		&\leq \int_{\Omega}|w|^{\alpha+\delta}dx + \frac{1}{2}\int_{\Omega}|w|^{2(\alpha + \delta - 1)}dx +\frac{1}{2}\int_{\Omega}|w|^2dx + C\int_{\Omega}z^2dx.
		\end{align}
		Since $\Phi(\nabla u_1,\nabla u_2) \geq 0$, \eqref{ine-u2} becomes
		\begin{align}\label{ine-u3}
		\frac{\lambda}{\alpha}\frac{d}{dt}\int_{\Omega}|w|^\alpha dx +\frac{1}{2}\frac{d}{dt}\int_{\Omega}z^2dx &\leq \int_{\Omega}|w|^{\alpha+\delta}dx+\frac{1}{2}\int_{\Omega}|w|^{2(\alpha + \delta - 1)}dx \nonumber\\&+ C\int_{\Omega}|w|^2dx + C\int_{\Omega}|z|^2dx.
		\end{align}
		We set $\delta=0$ and use the inequality \eqref{pre-ine-4} to rewrite \eqref{ine-u3} as
		\begin{align}\label{ine-u4}
		\frac{d}{dt}\left(\int_{\Omega}|w|^\alpha dx + \int_{\Omega}z^2dx\right) \leq C(\alpha,\lambda) +  C(\alpha,\lambda)\left(\int_{\Omega}|w|^\alpha dx + \int_{\Omega}z^2dx\right).
		\end{align}
		It is convenient to introduce the notation:
		$$W(t):=\int_{\Omega}|w|^\alpha dx + \int_{\Omega}|z|^2dx \text{ for } t \in S.$$
		Hence, the inequality \eqref{ine-u4} becomes
		\begin{align}\label{bfgronw}
		\frac{d}{dt}W(t) \leq C(\alpha,\lambda)W(t), 
		\end{align}
		for $t \in S$ with $W(0)=\int_{\Omega}|w_0|^\alpha dx + \int_{\Omega}|z_0|^2dx$, where $w_0:= u_{01}-u_{02}$ and $z_0:= v_{01}-v_{02}$. Here we consider $u_{01}, u_{02} \in L^\alpha(\Omega)$ and $v_{01}, v_{02} \in L^2(\Omega)$.\\	
		By using Gronwall's inequality, \eqref{bfgronw} yields
		\begin{align}
		W(t) \leq W(0)e^{tC(\alpha,\lambda)} \text{ for all } t \in S.
		\end{align}
		This also implies
		\begin{align}
		\int_{\Omega}|w|^\alpha dx + \int_{\Omega}|z|^2dx \leq (\|w_0\|_{L^{\alpha}(\Omega)}^\alpha + \|z_0\|_{L^2(\Omega)}^2)e^{tC(\alpha,\lambda)}.
		\end{align}
		Clearly, if $w_0=z_0=0$, then the weak solution of \eqref{main_eq} is unique.
	\end{proof}
	\section{Proof of Theorem \ref{stability-theo}}\label{stability}
	\begin{proof}
		Let us recall the weak formulation corresponding to the different choices of data: $(u_{0i}, v_{0i}, D_i, \varphi_i), i\in \{1,2\}$. We denote $D=D_1-D_2$, $\tilde{\varphi}= \varphi_1-\varphi_2$, $\tilde{r}=r_1-r_2$, $\tilde{u}_0 = u_{01}-u_{02}$ and $\tilde{v}_0 = v_{01}-v_{02}$. We denote also $w:=u_1 - u_2$ and  $z:= v_1-v_2$. Multiplying the first and the second equations of \eqref{main_eq} with $\psi := |w|^{\alpha + \delta - 1}, \phi := z$, respectively and interating the result by parts over $\Omega$ together with combining the two equations, one gets
		\begin{align}\label{stabi_1}
		&\int_{\Omega}\partial_t(u_1^\lambda - u_2^\lambda)\psi dx + \int_{\Omega}\partial_t(v_1-v_2)\phi dx + \int_{\Omega} \Big(K_1(|\nabla u_1|)\nabla u_1 \nonumber\\&- K_1(|\nabla u_2|)\nabla u_2\Big)\nabla \psi dx
		+\int_{\Omega}\left(D_1\nabla v_1 - D_2\nabla v_2\right)\nabla \phi dx + \int_{\Gamma^R}\left(\varphi_1u_1^\lambda - \varphi_2u_2^\lambda\right)\psi d\gamma \nonumber\\&=- \int_{\Omega}\left[r_1B(u_1 - v_1)(\psi - \phi) - r_2B(u_2-v_2)(\psi-\phi)\right]dx,
		\end{align}	
		Regarding \eqref{stabi_1}, note that 
		\begin{align}\label{pre-1}
		\int_{\Omega}\left(D_1\nabla v_1 - D_2\nabla v_2\right)\nabla \phi dx &= D_1\|\nabla \phi\|_{L^2(\Omega)}^2 + (D_1-D_2)\int_{\Omega}\nabla v_2 \nabla\phi dx,
		\end{align}	
		\begin{align}\label{pre-2}
		\int_{\Gamma^R}\left(\varphi_1 u_1^\lambda - \varphi_2u_2^\lambda\right)\psi d\gamma &= \int_{\Gamma^R} \varphi_1(u_1^\lambda - u_2^\lambda)\psi d\gamma + \int_{\Gamma^R}(\varphi_1-\varphi_2)u_2^\lambda\psi d\gamma
		\end{align}
		and
		\begin{align}\label{pre-3}
		&\int_{\Omega}\left[r_1B(u_1 - v_1)(\psi - \phi) - r_2B(u_2-v_2)(\psi-\phi)\right]dx \nonumber\\&= \int_{\Omega}r_1\left(B(u_1-v_1) - B(u_2-v_2)\right)(\psi - \phi)dx + (r_1-r_2)\int_{\Omega} B(u_2-v_2)(\psi - \phi)dx.
		\end{align}
		Using now \eqref{pre-1}, \eqref{pre-3}, as well as Young's inequality applied to the last terms of \eqref{pre-1}, \eqref{pre-2} together with the assumption $(\text{A}_2)$ and $(\text{A}_3)$, we use that \eqref{stabi_1} becomes
		\begin{align}\label{stabi_2}
		&\frac{\lambda}{\alpha}\frac{d}{dt}\int_{\Omega}|u_1-u_2|^\alpha dx + \frac{1}{2}\frac{d}{dt}\int_{\Omega}\phi^2dx +\int_{\Omega}\Phi(\nabla u_1,\nabla u_2)|u_1-u_2|^{\alpha+\delta-2}dx +\nonumber\\&+  D_1\int_{\Omega} |\nabla \phi|^2 dx + \int_{\Gamma^R} \varphi_1|u_1^\lambda - u_2^\lambda|\psi d\gamma 
		\leq C\varepsilon_1 \int_{\Omega}|\nabla \phi|^2dx  +\nonumber\\&+ C|D_1-D_2|\int_{\Omega}|\nabla v_2|^2dx - \frac{1}{2}\int_{\Gamma^R}|\varphi_1-\varphi_2|^2|u_2|^{2\lambda}d\gamma \nonumber\\&-\frac{1}{2}\int_{\Gamma^R} |u_1-u_2|^{2(\alpha+\delta-1)}d\gamma + \int_{\Omega}r_1(|u_1-u_2|+|v_1-v_2|)(\psi-\phi)dx  \nonumber\\&+ |r_1-r_2|\frac{\hat{c}}{\bar{r}}\int_{\Omega}|u_2-v_2|^\sigma (\psi-\phi)dx.
		\end{align} 
		Using the inequality \eqref{pre-ine-2}, \eqref{stabi_2} receives the form
		\begin{align}
		&\frac{\lambda}{\alpha}\frac{d}{dt}\int_{\Omega}|u_1-u_2|^\alpha dx + \frac{1}{2}\frac{d}{dt}\int_{\Omega}\phi^2dx +\int_{\Omega}\Phi(\nabla u_1,\nabla u_2)|u_1-u_2|^{\alpha+\delta-2}dx +\nonumber\\&+  D_1\int_{\Omega} |\nabla \phi|^2 dx + \int_{\Gamma^R} \varphi_1|u_1^\lambda - u_2^\lambda|\psi d\gamma 
		\leq C\varepsilon_1 \int_{\Omega}|\nabla \phi|^2dx  +\nonumber\\&+ C|D_1-D_2|\int_{\Omega}|\nabla v_2|^2dx - \frac{1}{2}\int_{\Gamma^R}|\varphi_1-\varphi_2|^2|u_2|^{2\lambda}d\gamma \nonumber\\&-\frac{1}{2}\int_{\Gamma^R} |u_1-u_2|^{2(\alpha+\delta-1)}d\gamma + \int_{\Omega}r_1|u_1-u_2|\psi dx - \int_{\Omega}r_1|u_1-u_2|\phi dx+ \nonumber\\& \int_{\Omega}r_1|v_1-v_2|\psi dx - \int_{\Omega}r_1|v_1-v_2|\psi dx + |r_1-r_2|\frac{\hat{c}}{\bar{r}}\int_{\Omega}(|u_2|^\sigma + |v_2|^\sigma)(\psi-\phi)dx.
		\end{align} 
		Applying the trace inequality \eqref{trace-theo} together with Cauchy-Schwarz's inequality, we obtain the following estimate
		\begin{align}
		&\frac{\lambda}{\alpha}\frac{d}{dt}\int_{\Omega}|u_1-u_2|^\alpha dx + \frac{1}{2}\frac{d}{dt}\int_{\Omega}\phi^2dx +\int_{\Omega}\Phi(\nabla u_1,\nabla u_2)|u_1-u_2|^{\alpha+\delta-2}dx \nonumber\\&+  D_1\int_{\Omega} |\nabla \phi|^2 dx + \int_{\Gamma^R} \varphi_1|u_1^\lambda - u_2^\lambda|\psi d\gamma 
		\leq C\varepsilon_1 \int_{\Omega}|\nabla \phi|^2dx  +\nonumber\\&+ C|D_1-D_2|\int_{\Omega}|\nabla v_2|^2dx - \frac{1}{2}\int_{\Gamma^R}|\varphi_1-\varphi_2|^2|u_2|^{2\lambda}d\gamma \nonumber\\&-\frac{1}{2}\Big(2\varepsilon_2 \int_{\Omega}|\nabla u_1 - \nabla u_2|^{2-a}|u_1-u_2|^{\alpha + \delta -2}dx + C\int_{\Omega}|u_1-u_2|^{\alpha}dx \nonumber\\&+ C\int_{\Omega}|u_1-u_2|^{\alpha+\mu_0}dx + C\int_{\Omega}|u_1-u_2|^{\alpha+\mu_1}dx\Big)\nonumber\\& +r_1\int_{\Omega}|u_1-u_2|^{\alpha+\delta}dx - \frac{r_1}{2}\int_{\Omega}|u_1-u_2|^2dx - \frac{r_1}{2}\int_{\Omega}|v_1-v_2|^2dx \nonumber\\& +\frac{r_1}{2}\int_{\Omega}|v_1-v_2|^2dx + \frac{r_1}{2}\int_{\Omega} |u_1-u_2|^{2(\alpha+\delta -1)}dx - r_1\int_{\Omega}|v_1-v_2|^2dx \nonumber\\& +|r_1-r_2|\frac{\hat{c}}{2\bar{r}}\int_{\Omega}|u_2|^{2\sigma}dx + |r_1-r_2|\frac{\hat{c}}{2\bar{r}}\int_{\Omega}|u_1-u_2|^{2(\alpha+\delta-1)}dx\nonumber\\&-|r_1-r_2|\frac{\hat{c}}{2\bar{r}}\int_{\Omega}|u_2|^{2\sigma}dx - |r_1-r_2|\frac{\hat{c}}{2\bar{r}}\int_{\Omega}|v_1-v_2|^{2}dx \nonumber\\& + |r_1-r_2|\frac{\hat{c}}{2\bar{r}}\int_{\Omega}|v_2|^{2\sigma}dx + |r_1-r_2|\frac{\hat{c}}{2\bar{r}}\int_{\Omega}|u_1-u_2|^{2(\alpha+\delta-1)}dx\nonumber\\&-|r_1-r_2|\frac{\hat{c}}{2\bar{r}}\int_{\Omega}|v_2|^{2\sigma}dx  - |r_1-r_2|\frac{\hat{c}}{2\bar{r}}\int_{\Omega}|v_1-v_2|^{2}dx.
		\end{align} 
		Choosing $\varepsilon_1=\frac{D_1}{C}$ and $\varepsilon_2=1$, we have 
		\begin{align}
		\frac{\lambda}{\alpha}\frac{d}{dt}\int_{\Omega}|u_1-u_2|^\alpha dx + \frac{1}{2}\frac{d}{dt}\int_{\Omega}\phi^2dx \leq C|D_1-D_2|\int_{\Omega}|\nabla v_2|^2dx -\nonumber\\ \frac{1}{2}\int_{\Gamma^R}|\varphi_1-\varphi_2|^2|u_2|^{2\lambda}d\gamma +r_1\int_{\Omega}|u_1-u_2|^{\alpha+\delta}dx - \frac{r_1}{2}\int_{\Omega}|u_1-u_2|^2dx \nonumber\\ + \frac{r_1}{2}\int_{\Omega} |u_1-u_2|^{2(\alpha+\delta -1)}dx - r_1\int_{\Omega}|v_1-v_2|^2dx \nonumber\\ + |r_1-r_2|\frac{\hat{c}}{\bar{r}}\int_{\Omega}|u_1-u_2|^{2(\alpha+\delta-1)}dx - |r_1-r_2|\frac{\hat{c}}{\bar{r}}\int_{\Omega}|v_1-v_2|^{2}dx.
		\end{align}
		Moreover, if we assume that $\delta = 0$, then the maximum allowed power of $\|w\|$ is $\alpha$. As next step, we use the inequality \eqref{pre-ine-5} together with the energy estimates \eqref{ine_main1}, \eqref{ine-main3} to deal with the terms $\int_{\Omega}|u_2|^\alpha dx$, $\int_{\Omega}|v_2|^2dx$ and $\int_{\Omega}|\nabla u_2|^2dx$. Furthermore, we use also the trace inequality \eqref{trace-theo}. It yields
		\begin{align}\label{stabi-5}
		\frac{d}{dt}\Big(\int_{\Omega}|u_1-u_2|^\alpha dx + \int_{\Omega}|v_1-v_2|^2dx\Big) \leq C(\alpha, \lambda, \hat{c}, \bar{r})(|D_1-D_2| + |r_1-r_2| \nonumber\\- \|\varphi_1-\varphi_2\|_{L^\infty}^2)+ C(\alpha, \lambda, \hat{c}, \bar{r})|r_1-r_2|\left(|u_1-u_2\|_{L^\alpha(\Omega)}^{\alpha} + \|v_1-v_2\|_{L^2(\Omega)}^{2}\right).
		\end{align}
		Denoting
		\begin{align}\label{stabi-6}
		Z(t):=\int_{\Omega}|u_1-u_2|^\alpha dx + \int_{\Omega}|v_1-v_2|^2dx \text{ for any } t \in S,
		\end{align}
		The expansion \eqref{stabi-5} can be rewritten as follows
		\begin{align}
		\frac{d}{dt}Z(t) \leq C(\alpha, \lambda, \hat{c}, \bar{r})(|D_1-D_2| + |r_1-r_2| - \|\varphi_1-\varphi_2\|_{L^\infty(\Gamma^R)}^2) \nonumber\\+ C(\alpha, \lambda, \hat{c}, \bar{r})|r_1-r_2|Z(t),
		\end{align}
		for $t\in S$. It holds $Z(0) = \int_{\Omega}|u_{01}-u_{02}|^\alpha dx + \int_{\Omega}|v_{01}-v_{02}|^2dx$.\\
		Appying the Gr\"{o}nwall's inequality to \eqref{stabi-6}, we obtain
		\begin{align}\label{stabi-7}
		Z(t) \leq e^{\int_{0}^tC(\alpha,\lambda,\hat{c},\bar{r})|r_1-r_2|ds}\Big[Z(0) + 
		\int_0^tC(\alpha, \lambda, \hat{c}, \bar{r})(|D_1-D_2|+ |r_1-r_2| \nonumber\\- \|\varphi_1-\varphi_2\|_{L^\infty(\Gamma^R)}^2 )ds\Big].
		\end{align}
		\eqref{stabi-7} implies
		\begin{align}
		\|u_1-u_2\|_{L^\alpha(\Omega)}^\alpha + \|v_1-v_2\|_{L^2(\Omega)}^2 \leq e^{C(\alpha,\lambda,\hat{c},\bar{r})|r_1-r_2|t}\Big[ \|u_{01} - u_{02}\|_{L^\alpha(\Omega)}^\alpha \nonumber\\+ \|v_{01} - v_{02}\|_{L^2(\Omega)}^2 + Ct(|D_1-D_2|+ |r_1-r_2| - \|\varphi_1-\varphi_2\|_{L^\infty(\Gamma^R)}^2 )\Big],
		\end{align}
		which is precisely the kind of stability estimate with respect to data and parameters we are looking for.
	\end{proof}
	\section{Proof of Theorem \ref{stabilityfull-theo}}\label{stability-full}
	\begin{proof}
		We keep the same notations as in Theorem \ref{stability-theo}. Multiplying the first and the second equation of \eqref{main_eq} with $\psi:= (u_1-u_2)_t$ and $\phi:=(v_1-v_2)_t$, respectively. Integrating the results by parts over $\Omega$ and combining the two equations, we obtain
		\begin{align}\label{stabi2-1}
		&\int_{\Omega}\partial_t(u_1^\lambda - u_2^\lambda)\psi dx + \int_{\Omega}\partial_t(v_1-v_2)\phi dx + \int_{\Omega} \Big(K_1(|\nabla u_1|)\nabla u_1 \nonumber\\&- K_1(|\nabla u_2|)\nabla u_2\Big)\nabla \psi dx
		+\int_{\Omega}\left(D_1\nabla v_1 - D_2\nabla v_2\right)\nabla \phi dx  \nonumber\\&+ \int_{\Gamma^R}\varphi\left( u_1^\lambda - u_2^\lambda\right)\psi d\gamma =- \int_{\Omega}\Big[r_1B(u_1 - v_1)(\psi - \phi) \nonumber\\&- r_2B(u_2-v_2)(\psi-\phi)\Big]dx.
		\end{align}
		Note that 
		\begin{align}\label{stabi-full1}
		\int_{\Omega}(D_1\nabla v_1 - D_2\nabla v_2)\nabla \phi dx = \int_{\Omega}D_1 \nabla (v_1-v_2)\nabla \phi dx + \int_{\Omega}(D_1-D_2)\nabla v_2 \nabla \phi dx \nonumber\\
		= \frac{1}{2}\frac{d}{dt}\int_{\Omega} D_1|\nabla(v_1-v_2)|^2dx + \int_{\Omega}(D_1-D_2)\nabla v_2 \nabla \phi dx,
		\end{align}
		\begin{align}\label{stabi-full2}
		\int_{\Gamma^R} (\varphi_1u_1^\lambda - \varphi_2 u_2^\lambda)\psi d\gamma = \int_{\Gamma^R} \varphi_1(u_1^\lambda - u_2^\lambda)\psi d\gamma + \int_{\Gamma^R} (\varphi_1-\varphi_2)u_2^\lambda \psi d\gamma
		\end{align}
		and
		\begin{align}\label{stabi-full3}
		&\int_{\Omega}\left[r_1B(u_1- v_1)(\psi - \phi) - r_2B(u_2-v_2)(\psi - \phi)\right]dx = \nonumber\\ &\int_{\Omega}r_1(B(u_1-v_1) - B(u_2-v_2))(\psi - \phi)dx + (r_1-r_2)\int_{\Omega}B(u_2-v_2)(\psi-\phi)dx.
		\end{align}
		Then, \eqref{stabi2-1} becomes 
		\begin{align}
		\frac{d}{dt}\int_{\Omega}(u_1-u_2)^2 dx + \int_{\Omega}(K_1(|\nabla u_1|)\nabla u_1 - K_1(|\nabla u_2|)\nabla u_2)\cdot \frac{\partial}{\partial t}(\nabla u_1 - \nabla u_2) dx \nonumber\\
		+ \frac{1}{2}\frac{d}{dt}\int_{\Omega}|v_1-v_2|^2dx + \int_{\Omega}(D_1\nabla v_1-D_2\nabla v_2)\nabla \phi dx + \int_{\Gamma^R}(\varphi_1u_1^\lambda - \varphi_2u_2^\lambda)\psi d\gamma \nonumber\\= \int_{\Omega}\left[r_1B(u_1- v_1)(\psi - \phi) - r_2B(u_2-v_2)(\psi - \phi)\right]dx\nonumber\\+ (r_1-r_2)\int_{\Omega}B(u_2-v_2)(\psi-\phi)dx.
		\end{align}
		Using \eqref{stabi-full1}-\eqref{stabi-full3}, as a result of applying the Young's inequality to the last terms of the right hand side of \eqref{stabi-full1} and \eqref{stabi-full2}, we have
		\begin{align}\label{stabi-full4}
		&\frac{d}{dt}\int_{\Omega}(u_1-u_2)^2 dx + \int_{\Omega}K_1(|\nabla u_1|)\nabla u_1\cdot \frac{\partial }{\partial t}(\nabla u_1)dx + \int_{\Omega}K_1(|\nabla u_2|)\nabla u_2\cdot \frac{\partial }{\partial t}(\nabla u_2)dx \nonumber\\&- \int_{\Omega}K_1(|\nabla u_1|)\nabla u_1\cdot \frac{\partial }{\partial t}(\nabla u_2)dx - \int_{\Omega}K_1(|\nabla u_2|)\nabla u_2\cdot \frac{\partial }{\partial t}(\nabla u_1)dx \nonumber\\&+ \frac{1}{2}\frac{d}{dt}\int_{\Omega}|v_1-v_2|^2dx+ \frac{1}{2}\frac{d}{dt}\int_{\Omega} D_1|\nabla (v_1-v_2)|^2dx + \int_{\Gamma^R}\varphi_1|u_1^\lambda - u_2^\lambda|\psi d\gamma \nonumber\\
		&\leq C\varepsilon_3 \frac{d}{dt}\int_{\Omega}|\nabla v_1 - \nabla v_2|^2dx + C|D_1-D_2|\int_{\Omega}|\nabla v_2|^2dx   \nonumber\\&- \frac{1}{2}\|\varphi_1-\varphi_2\|_{L^\infty(\Gamma^R)}^2\int_{\Gamma^R}|u_2|^{2\lambda}d\gamma - \frac{1}{2}\frac{d}{dt}\int_{\Gamma^R}|u_1-u_2|^2d\gamma + \int_{\Omega}r_1B(u_1-v_1)\psi dx  \nonumber\\&- \int_{\Omega}r_1B(u_1-v_1)\phi dx - \int_{\Omega}r_2B(u_2-v_2)\psi dx  \nonumber\\&+ \int_{\Omega}r_2B(u_2-v_2)\phi dx  + |r_1-r_2|\int_{\Omega}B(u_2-v_2)\psi dx  \nonumber\\&- |r_1-r_2|\int_{\Omega}B(u_2-v_2)\phi dx.
		\end{align}
		By the assumption $(\text{A}_3)$, we are led to
		\begin{align}\label{stabi-full6}
		&\frac{d}{dt}\int_{\Omega}(u_1-u_2)^2 dx + \frac{d}{dt}\int_{\Omega}K_1(|\nabla u_1|)|\nabla u_1|^2dx + \frac{d}{dt}\int_{\Omega}K_1(|\nabla u_2|)|\nabla u_2|^2dx \nonumber\\&- \frac{d}{dt}\int_{\Omega}K_1(|\nabla u_1|)\nabla u_1\cdot \nabla u_2dx - \frac{d}{dt}\int_{\Omega}K_1(|\nabla u_2|)\nabla u_2\cdot\nabla u_1dx \nonumber\\&+ \frac{1}{2}\frac{d}{dt}\int_{\Omega}|v_1-v_2|^2dx+ \frac{1}{2}\frac{d}{dt}\int_{\Omega} D_1|\nabla (v_1-v_2)|^2dx + \int_{\Gamma^R}\varphi_1|u_1^\lambda - u_2^\lambda|\psi d\gamma \nonumber\\
		&\leq C\varepsilon_3 \frac{d}{dt}\int_{\Omega}|\nabla v_1 - \nabla v_2|^2dx + C|D_1-D_2|\int_{\Omega}|\nabla v_2|^2dx  \nonumber\\&- \frac{1}{2}\|\varphi_1-\varphi_2\|_{L^\infty(\Gamma^R)}^2\int_{\Gamma^R}|u_2|^{2\lambda}d\gamma - \frac{1}{2}\frac{d}{dt}\int_{\Gamma^R}|u_1-u_2|^2d\gamma  \nonumber\\&+ \int_{\Omega} r_1(B(u_1-v_1) - B(u_2-v_2))\psi dx + \int_{\Omega}r_1(B(u_1-v_1) - B(u_2-v_2))\phi dx  \nonumber\\&+ \frac{\hat{c}}{\bar{r}}|r_1-r_2|\int_{\Omega}|u_2-v_2|^\sigma \psi dx - \frac{\hat{c}}{\bar{r}}|r_1-r_2|\int_{\Omega}|u_2-v_2|^\sigma \phi dx
		\end{align}
		By using the assumption $(\text{A}_2)$ together with \eqref{pre-ine-2}, it leads to
		\begin{align}\label{stabi-full8}
		&\frac{d}{dt}\int_{\Omega}(u_1-u_2)^2 dx + \frac{d}{dt}\int_{\Omega}K_1(|\nabla u_1|)|\nabla u_1|^2dx + \frac{d}{dt}\int_{\Omega}K_1(|\nabla u_2|)|\nabla u_2|^2dx \nonumber\\&- \frac{d}{dt}\int_{\Omega}K_1(|\nabla u_1|)\nabla u_1\cdot \nabla u_2dx - \frac{d}{dt}\int_{\Omega}K_1(|\nabla u_2|)\nabla u_2\cdot\nabla u_1dx \nonumber\\&+ \frac{1}{2}\frac{d}{dt}\int_{\Omega}|v_1-v_2|^2dx+ \frac{1}{2}\frac{d}{dt}\int_{\Omega} D_1|\nabla (v_1-v_2)|^2dx + \int_{\Gamma^R}\varphi_1|u_1^\lambda - u_2^\lambda|\psi d\gamma \nonumber\\
		&\leq C\varepsilon_3 \frac{d}{dt}\int_{\Omega}|\nabla v_1 - \nabla v_2|^2dx + C|D_1-D_2|\int_{\Omega}|\nabla v_2|^2dx  \nonumber\\& - \frac{1}{2}\|\varphi_1-\varphi_2\|_{L^\infty(\Gamma^R)}^2\int_{\Gamma^R}|u_2|^{2\lambda}d\gamma - \frac{1}{2}\frac{d}{dt}\int_{\Gamma^R}|u_1-u_2|^2d\gamma \nonumber\\&+ \int_{\Omega}r_1|u_1-v_1-u_2+v_2|\psi dx + \int_{\Omega}r_1|u_1-v_1-u_2+v_2|\phi dx +\frac{\hat{c}}{\bar{r}}|r_1-r_2|\int_{\Omega}(|u_2|^\sigma \nonumber\\&+ |v_2|^\sigma)\psi dx - \frac{\hat{c}}{\bar{r}}|r_1-r_2|\int_{\Omega}(|u_2|^\sigma + |v_2|^\sigma)\phi dx.
		\end{align}
		Then, via Cauchy-Schwarz's inequality, we obtain the following estimate
		\begin{align}\label{stabi-full9}
		&\frac{d}{dt}\int_{\Omega}(u_1-u_2)^2 dx+ \frac{d}{dt}\int_{\Omega}K_1(|\nabla u_1|)|\nabla u_1|^2dx + \frac{d}{dt}\int_{\Omega}K_1(|\nabla u_2|)|\nabla u_2|^2dx \nonumber\\&- \frac{d}{dt}\int_{\Omega}K_1(|\nabla u_1|)\nabla u_1\cdot \nabla u_2dx - \frac{d}{dt}\int_{\Omega}K_1(|\nabla u_2|)\nabla u_2\cdot\nabla u_1dx  \nonumber\\&+ \frac{1}{2}\frac{d}{dt}\int_{\Omega}|v_1-v_2|^2dx+ \frac{1}{2}\frac{d}{dt}\int_{\Omega} D_1|\nabla (v_1-v_2)|^2dx + \int_{\Gamma^R}\varphi_1|u_1^\lambda - u_2^\lambda|\psi d\gamma  \nonumber\\
		&\leq C\varepsilon_3 \frac{d}{dt}\int_{\Omega}|\nabla v_1 - \nabla v_2|^2dx+ C|D_1-D_2|\int_{\Omega}|\nabla v_2|^2dx  \nonumber\\& - \frac{1}{2}\|\varphi_1-\varphi_2\|_{L^\infty(\Gamma^R)}^2\int_{\Gamma^R}|u_2|^{2\lambda}d\gamma - \frac{1}{2}\frac{d}{dt}\int_{\Gamma^R}|u_1-u_2|^2d\gamma \nonumber\\&+ r_1\int_{\Omega}|u_1-u_2|\psi dx + r_1\int_{\Omega}|v_1-v_2|\psi dx + r_1\int_{\Omega}|u_1-u_2|\phi dx + r_1\int_{\Omega}|v_1-v_2|\phi dx \nonumber\\&+ \frac{\hat{c}}{2\bar{r}}|r_1-r_2|\int_{\Omega}|u_2|^{2\sigma}dx  + \frac{\hat{c}}{4\bar{r}}|r_1-r_2|\frac{d}{dt}\int_{\Omega}|u_1-u_2|^2dx \nonumber\\&+ \frac{\hat{c}}{2\bar{r}}|r_1-r_2|\int_{\Omega}|v_2|^{2\sigma}dx +\frac{\hat{c}}{4\bar{r}}|r_1-r_2|\frac{d}{dt}\int_{\Omega}|u_1-u_2|^2dx \nonumber\\&- \frac{\hat{c}}{2\bar{r}}|r_1-r_2|\int_{\Omega}|u_2|^{2\sigma}dx - \frac{\hat{c}}{4\bar{r}}|r_1-r_2|\frac{d}{dt}\int_{\Omega}|v_1-v_2|^2dx  \nonumber\\&- \frac{\hat{c}}{2\bar{r}}|r_1-r_2|\int_{\Omega}|v_2|^{2\sigma}dx - \frac{\hat{c}}{4\bar{r}}|r_1-r_2|\frac{d}{dt}\int_{\Omega}|v_1-v_2|^2dx.
		\end{align}
		In other words, \eqref{stabi-full8} can be written as follow
		\begin{align}\label{stabi-full10}
		&\frac{d}{dt}\int_{\Omega}(u_1-u_2)^2 dx+ \frac{d}{dt}\int_{\Omega}K_1(|\nabla u_1|)|\nabla u_1|^2dx + \frac{d}{dt}\int_{\Omega}K_1(|\nabla u_2|)|\nabla u_2|^2dx \nonumber\\&- \frac{d}{dt}\int_{\Omega}K_1(|\nabla u_1|)\nabla u_1\cdot \nabla u_2dx - \frac{d}{dt}\int_{\Omega}K_1(|\nabla u_2|)\nabla u_2\cdot\nabla u_1dx  \nonumber\\&+ \frac{1}{2}\frac{d}{dt}\int_{\Omega}|v_1-v_2|^2dx+ \frac{1}{2}\frac{d}{dt}\int_{\Omega} D_1|\nabla (v_1-v_2)|^2dx + \int_{\Gamma^R}\varphi_1|u_1^\lambda - u_2^\lambda|\psi d\gamma \nonumber\\
		&\leq C\varepsilon_3 \frac{d}{dt}\int_{\Omega}|\nabla v_1 - \nabla v_2|^2dx + C|D_1-D_2|\int_{\Omega}|\nabla v_2|^2dx   \nonumber\\& - \frac{1}{2}\|\varphi_1-\varphi_2\|_{L^\infty(\Gamma^R)}^2\int_{\Gamma^R}|u_2|^{2\lambda}d\gamma - \frac{1}{2}\frac{d}{dt}\int_{\Gamma^R}|u_1-u_2|^2d\gamma + r_1\int_{\Omega}|u_1-u_2|^2dx \nonumber\\& + \frac{1}{2}\frac{d}{dt}\int_{\Omega}|u_1-u_2|^2dx + r_1\int_{\Omega}|v_1-v_2|^2dx + \frac{r_1}{2}\frac{d}{dt}\int_{\Omega}|v_1-v_2|^2dx \nonumber\\&+ |r_1-r_2|\frac{\hat{c}}{2\bar{r}}\frac{d}{dt}\int_{\Omega}|u_1-u_2|^2dx - |r_1-r_2|\frac{\hat{c}}{2\bar{r}}\frac{d}{dt}\int_{\Omega}|v_1-v_2|^2dx.
		\end{align}
		Using the property of $H$ in \eqref{H-ine1}, \eqref{stabi-full10} becomes
		\begin{align}\label{stabi-full11}
		&\frac{d}{dt}\int_{\Omega}(u_1-u_2)^2 dx + \frac{1}{2}\frac{d}{dt}\int_{\Omega}(H(|\nabla u_1|) + H(|\nabla u_2|))dx + \frac{1}{2}\frac{d}{dt}\int_{\Omega}|v_1-v_2|^2dx \nonumber\\&  + \frac{D_1}{2}\frac{d}{dt}\int_{\Omega}|\nabla u_1-\nabla u_2|^2dx \leq  \frac{d}{dt}\int_{\Omega}K_1(|\nabla u_1|)\nabla u_1\cdot \nabla u_2dx \nonumber\\&+ \frac{d}{dt}\int_{\Omega}K_1(|\nabla u_2|)\nabla u_2\cdot\nabla u_1dx  +C\varepsilon_3 \frac{d}{dt}\int_{\Omega}|\nabla v_1 - \nabla v_2|^2dx \nonumber\\&+ C|D_1-D_2|\int_{\Omega}|\nabla v_2|^2dx  - \frac{1}{2}\|\varphi_1-\varphi_2\|_{L^\infty(\Gamma^R)}^2\int_{\Gamma^R}|u_2|^{2\lambda}d\gamma \nonumber\\& - \frac{1}{2}\frac{d}{dt}\int_{\Gamma^R}|u_1-u_2|^2d\gamma + r_1\int_{\Omega}|u_1-u_2|^2dx + \frac{1}{2}\frac{d}{dt}\int_{\Omega}|u_1-u_2|^2dx \nonumber\\& + r_1\int_{\Omega}|v_1-v_2|^2dx + \frac{r_1}{2}\frac{d}{dt}\int_{\Omega}|v_1-v_2|^2dx + |r_1-r_2|\frac{\hat{c}}{2\bar{r}}\frac{d}{dt}\int_{\Omega}|u_1-u_2|^2dx \nonumber\\&- |r_1-r_2|\frac{\hat{c}}{2\bar{r}}\frac{d}{dt}\int_{\Omega}|v_1-v_2|^2dx.
		\end{align}
		Integrating \eqref{stabi-full11} over time interval $(0,t)$, it yields
		\begin{align}\label{stabi-full12}
		&\int_{\Omega}(u_1-u_2)^2 dx + \frac{1}{2}\int_{\Omega}(H(|\nabla u_1|) + H(|\nabla u_2|))dx \nonumber\\&+ \frac{D_1}{2}\int_{\Omega}|\nabla u_1-\nabla u_2|^2dx \leq \int_{\Omega}K_1(|\nabla u_1|)\nabla u_1\cdot \nabla u_2dx \nonumber\\&+ \int_{\Omega}K_1(|\nabla u_2|)\nabla u_2\cdot\nabla u_1dx + C\varepsilon_3\int_{\Omega}|\nabla v_1 - \nabla v_2|^2dx  \nonumber\\&+ C|D_1-D_2|\int_0^t\int_{\Omega}|\nabla v_2|^2dxds  - \frac{1}{2}\|\varphi_1-\varphi_2\|_{L^\infty(\Gamma^R)}^2\int_0^t\int_{\Gamma^R}|u_2|^{2\lambda}d\gamma ds \nonumber\\&- \frac{1}{2}\int_{\Gamma^R}|u_1-u_2|^2d\gamma  + r_1\int_0^t\int_{\Omega}|u_1-u_2|^2dxds + \frac{1}{2}\int_{\Omega}|u_1-u_2|^2dx \nonumber\\&+ r_1\int_0^t\int_{\Omega}|v_1-v_2|^2dxds  + \frac{r_1}{2}\int_{\Omega}|v_1-v_2|^2dx+ |r_1-r_2|\frac{\hat{c}}{2\bar{r}}\int_{\Omega}|u_1-u_2|^2dx\nonumber\\& - |r_1-r_2|\frac{\hat{c}}{2\bar{r}}\int_{\Omega}|v_1-v_2|^2dx.
		\end{align}
		Next, we estimate the first term of the right hand side of \eqref{stabi-full12} by using assumption $(\text{A}_5)$ together with Hold\"{e}r's inequality, one obtains
		\begin{align}\label{time-est1}
		\left|\int_{\Omega}K_1(|\nabla u_1|)\nabla u_1 \cdot \nabla u_2 dx\right| &\leq C\int_{\Omega}|\nabla u_1|^{1-a}|\nabla u_2|dx \nonumber\\
		&\leq C\left(\int_{\Omega}|\nabla u_1|^{2-a}dx\right)^{\frac{1-a}{2-a}}\left(\int_{\Omega}|\nabla u_2|^{2-a}dx\right)^{\frac{1}{2-a}} \nonumber\\
		&\leq C\left(\int_{\Omega}(H(|\nabla u_1|) + 1)dx\right)^{\frac{1-a}{2-a}}\left(\int_{\Omega}|\nabla u_2|^{2-a}dx\right)^{\frac{1}{2-a}} \nonumber\\
		&\leq C\left(\int_{\Omega}H(|\nabla u_1|)dx\right)^{\frac{1-a}{2-a}}\left(\int_{\Omega}|\nabla u_2|^{2-a}dx\right)^{\frac{1}{2-a}} \nonumber\\
		&+ C\left(\int_{\Omega}|\nabla u_2|^{2-a}dx\right)^{\frac{1}{2-a}} \nonumber\\
		&\leq \varepsilon \int_{\Omega}H(|\nabla u_1|)dx + C\int_{\Omega}|\nabla u_2|^{2-a}dx \nonumber\\
		&+ C\left(\int_{\Omega}|\nabla u_2|^{2-a}dx\right)^{\frac{1}{2-a}}.
		\end{align}
		By using the same procedure for the second term of the right hand side of \eqref{stabi-full12}, we have
		\begin{align}\label{time-est2}
		\left|\int_{\Omega}K_1(|\nabla u_2|)\nabla u_2 \cdot \nabla u_1 dx\right| \leq \varepsilon \int_{\Omega}H(|\nabla u_2|)dx + C\int_{\Omega}|\nabla u_1|^{2-a}dx \nonumber\\+ C\left(\int_{\Omega}|\nabla u_1|^{2-a}dx\right)^{\frac{1}{2-a}}.
		\end{align}
		Applying trace inequality \eqref{trace-theo} and using \eqref{time-est1}, \eqref{time-est2}, we obtain the following estimate
		\begin{align}
		&\int_{\Omega}(u_1-u_2)^2 dx + \frac{1}{2}\int_{\Omega}(H(|\nabla u_1|) + H(|\nabla u_2|))dx \nonumber\\& + \frac{D_1}{2}\int_{\Omega}|\nabla u_1-\nabla u_2|^2dx \leq \varepsilon \int_{\Omega}(H(|\nabla u_1|) + H(|\nabla u_2|))dx \nonumber\\&+ C\int_{\Omega}|\nabla u_1|^{2-a}dx + C\left(\int_{\Omega}|\nabla u_1|^{2-a}dx\right)^{\frac{1}{2-a}} \nonumber\\&+ C
		\int_{\Omega}|\nabla u_2|^{2-a}dx + C\left(\int_{\Omega}|\nabla u_2|^{2-a}dx\right)^{\frac{1}{2-a}} \nonumber\\&+ C\varepsilon_3\int_{\Omega}|\nabla v_1 - \nabla v_2|^2dx + C|D_1-D_2|\int_0^t\int_{\Omega}|\nabla v_2|^2dxds - \nonumber\\&  \frac{1}{2}\|\varphi_1-\varphi_2\|_{L^\infty(\Gamma^R)}^2\int_0^t\int_{\Gamma^R}|u_2|^{2\lambda}d\gamma ds - \frac{1}{2}\Big(2\varepsilon_4\int_{\Omega}|\nabla u_1-\nabla u_2|^{2-a}dx \nonumber\\&+ C\int_{\Omega}|u_1-u_2|^\alpha dx + C\int_{\Omega}|u_1-u_2|^{\alpha+\mu_0}dx + C\int_\Omega|u_1-u_2|^{\alpha+\mu_1}dx\Big) \nonumber\\&+ r_1\int_0^t\int_{\Omega}|u_1-u_2|^2dxds  + \frac{1}{2}\int_{\Omega}|u_1-u_2|^2dx + r_1\int_0^t\int_{\Omega}|v_1-v_2|^2dxds  \nonumber\\&+ \frac{r_1}{2}\int_{\Omega}|v_1-v_2|^2dx+ |r_1-r_2|\frac{\hat{c}}{2\bar{r}}\int_{\Omega}|u_1-u_2|^2dx  \nonumber\\&- |r_1-r_2|\frac{\hat{c}}{2\bar{r}}\int_{\Omega}|v_1-v_2|^2dx.
		\end{align}
		Choosing $\varepsilon=\frac{1}{2}$, $\varepsilon_3=\frac{1}{4C}$ and $\varepsilon_4=1$, we are led to
		\begin{align}
		& \int_{\Omega}|\nabla u_1-\nabla u_2|^{2-a}dx + \int_{\Omega}|\nabla v_1-\nabla v_2|^2dx \leq  \varepsilon \int_{\Omega}(H(|\nabla u_1|) + H(|\nabla u_2|))dx \nonumber\\&+ C\int_{\Omega}|\nabla u_1|^{2-a}dx + C\left(\int_{\Omega}|\nabla u_1|^{2-a}dx\right)^{\frac{1}{2-a}} \nonumber\\&+ C
		\int_{\Omega}|\nabla u_2|^{2-a}dx + C\left(\int_{\Omega}|\nabla u_2|^{2-a}dx\right)^{\frac{1}{2-a}} \nonumber\\&+C|D_1-D_2|\int_0^t\int_{\Omega}|\nabla v_2|^2dxds -\frac{1}{2}\|\varphi_1-\varphi_2\|_{L^\infty(\Gamma^R)}^2\int_0^t\int_{\Gamma^R}|u_2|^{2\lambda}d\gamma ds\nonumber\\& + r_1\int_0^t\int_{\Omega}|u_1-u_2|^2dxds + \frac{1}{2}\int_{\Omega}|u_1-u_2|^2dx + r_1\int_0^t\int_{\Omega}|v_1-v_2|^2dxds +\nonumber\\& \frac{r_1}{2}\int_{\Omega}|v_1-v_2|^2dx + |r_1-r_2|\frac{\hat{c}}{2\bar{r}}\int_{\Omega}|u_1-u_2|^2dx - |r_1-r_2|\frac{\hat{c}}{2\bar{r}}\int_{\Omega}|v_1-v_2|^2dx.
		\end{align}
		Now, we use the energy estimates \eqref{ine_main1}, \eqref{ine-main3} and the stability-like estimate \eqref{stabi-main1}, we finally obtain the following structural-stability-like estimate
		\begin{align}
		\int_{\Omega}|\nabla u_1-\nabla u_2|^{2-a}dx + \int_{\Omega}|\nabla v_1-\nabla v_2|^2dx \leq C+ C(|D_1-D_2| -\nonumber\\ \|\varphi_1-\varphi_2\|_{L^\infty(\Gamma^R)}^2)  + \left(\left(t + \frac{1}{2}\right) + |r_1-r_2|\frac{\hat{c}}{2\bar{r}}\right)e^{C(\alpha,\lambda,\hat{c},\bar{r})|r_1-r_2|t}\Big(\|u_{01} - u_{02}\|_{L^2(\Omega)}^2 \nonumber\\ + \|v_{01}-v_{02}\|_{L^2(\Omega)}^2+ Ct\left(|D_1-D_2|+|r_1-r_2| - \|\varphi_1-\varphi_2\|_{L^\infty(\Gamma^R)}^2\right)\Big).
		\end{align}
	\end{proof}
	
	
	\begin{center}
		\bibliographystyle{alpha}
		\bibliography{mybib}
	\end{center}

	\medskip
	\medskip
	
\end{document}